\documentclass{kms-b}

\usepackage{color,graphicx}
\usepackage{latexsym}
\usepackage{amsmath,amssymb,amsthm,amsfonts}
\usepackage[all]{xy}
\usepackage{array}

\newtheorem{lem}{Lemma}
\newtheorem{theorem}{Theorem}

\newtheorem{cor}{Corollary}
\newtheorem{definition}{Definition}

\newtheorem{thm}{Theorem}

\newtheorem{prop}{Proposition}

\newcommand{\Q}{{\mathbb Q}}
\newcommand{\C}{{\mathbb C}}
\newcommand{\s}{{\sigma}}
\newcommand{\e}{{\varepsilon}}
\newcommand{\Z}{{\mathbb Z}}

\begin{document}
\title[ Certain homotopy properties related to $\text{map}(\Sigma^n \C P^2,S^m)$]
{ Certain homotopy properties related to $\text{map}(\Sigma^n \C P^2,S^m)$}

\author{Jin-ho Lee}
\address{Jin Ho Lee\\
Department of Mathematics\\
Korea University\\
Seoul 136-701, Korea}
\email{sabforev@korea.ac.kr}


\begin{abstract}
We compute cohomotopy groups of suspended complex plane $\pi^{n+m}(\Sigma^n \C P^2)$
for $m=6, 7, 8$.
Using
these results, we classify path components of the spaces
$map(\Sigma^n \C P^2,S^m)$ up to homotopy equivalent.
We determine the generalized Gottlieb groups $G_n(\C P^2,S^m)$.
Finally, we compute homotopy groups of mapping spaces
$map(\Sigma^n \mathbb{C}P^2,S^m;f)$ for all generators $[f]$ of
$[\Sigma^n \C P^2,S^m]$ and Gottlieb groups of mapping components containing constant map
$map(\Sigma^n \C P^2,S^m;\ast)$.
\end{abstract}

\maketitle

\section{Introduction}

Let $X$ and $Y$ be based topological spaces. A major object of homotopy theory is to study $[X,Y]$,
the set of homotopy equivalence classes of based maps.
In general, if $Y$ is a co-H-group, $[X,Y]$ has group structure.
Let $\Sigma X$ be the suspension of $X$.
Since every suspended space $\Sigma X$ is co-group, $[\Sigma X,Y]$ has group structure.
If $\Sigma X$ is a sphere $S^n$, $[S^n,Y]$ is the $n$-th homotopy group of $Y$.
On the other hand, $[X,S^n]$ is called the $n$-th cohomotopy set of $X$ and denoted by $\pi^n(X)$.
Moreover, if $X$ is a co-H-group, the cohomotopy set is a group and called $n$-th cohomotopy group of $X$.
Homotopy groups and cohomotopy groups have been studied by many authors and are a major object in algebraic topology.

Another major object of homotopy theory is to investigate the set of (unbased) maps $f:X \to Y$.
%
We denote $map(X,Y)$ to be the set of all continuous maps from $X$ to $Y$ equipped with compact-open topology.
Then we write $map(X,Y;f)$ for the path-component of $map(X,Y)$ containing a map $f:X \to Y$.
Important cases are $map(X,Y;\ast)$ and $map(X,X;1)$ where $\ast :X \to Y$ is the constant map and $1:X \to X$ is the identity map on $X$.
D. Harris proved that every topological space $X$ is the space of path components of a stratifiable space $S(X)$ \cite{HA}.
Thus it is important to study path components of mapping spaces
to analyze topological spaces.

Lang proved that if $X$ is a suspended space, then all path
components of based mapping space $map_\ast(X,Y)$ have the same homotopy type
\cite[Theorem 2.1]{Lang}.
Whitehead proved that $map(S^n,S^m;f)$
is homotopy equivalent to $map(S^n,S^m;\ast)$ if and only if $w_f$
has a section, where $\ast: S^n \to S^m$ is a constant map \cite[Theorem 2.8]{White}.
Smith proved a rational homotopy equivalence
$$map(X,X;\ast) \simeq_\Q X \times map_\ast(X,X:\ast)$$
for a ``nice'' space $X$ \cite{Sm}.
Lupton and Smith proved that
$$map(X,Y;f) \simeq map(X,Y;f+d)$$
for a CW co-H-space $X$ and any
CW complex $Y$, where $d:X \to Y$ is a cyclic map \cite[Theorem 3.10]{LS}.
%
Recently, Gatsinzi \cite{GAT} proved that the dimension of the rational
Gottlieb group of the universal cover $\widetilde{map}(X,S^{2n};f)$ of the function space $map(X,S^{2n};f)$ is at least equal to the dimension of $\tilde{H}^\ast(X;\Q)$ under several assumptions.
Lupton and Smith \cite{LS} showed that
$$G_n(map_\ast(X,Y;\ast)) \cong G_n(Y) \oplus G_n(X,Y).$$
Maruyama and Oshima \cite{MO} determined homotopy groups of $map_\ast(SU(3),SU(3);\ast)$, $map_\ast(Sp(2),Sp(2);\ast)$ and $map_\ast(G_2,G_2;\ast)$.
\\

In Section 2, we present some basic knowledge of composition methods \cite{T}.
We review a mapping cone sequence and a Puppe sequence related
to the suspended complex projective plane and discuss the concept
of cyclic maps and its properties.
Also we recall the Toda bracket and its properties.
\\

In Sections 3, 4 and 5, we compute $\pi^n(\Sigma^{n+k} \C P^2)$
for $k=6, 7$. As a result, we obtain the following results:

 {\begin{table}[htbp]
 \setlength\extrarowheight{2.5pt}
 \centering
\begin{tabular}{|c||c|c|c|c|c|c|c|c|c|c|}
\hline
$case~n$ &2 & 3  & 4& 5  & 6 &7 \\
\hline
$~$ &&&&&\\[-3.7ex]
$\pi^n(\Sigma^{n+6}\C P^2)$ & $2+15$ & $2+3 $& $8+2+3^2+5$ & $4+9 $ & $4^2+9+3 $ & $4+3 $ \\
\hline
$case~n$  &8 &9 &10 &11  & $n\geq 12$ & \\
\hline  $~$ &&&&&\\[-3.3ex]
$\pi^n(\Sigma^{n+6}\C P^2)$   &$ 4^2+3^2$ &$ 4+3 $& $2+3 $& $ 2+3$ & 2 &\\
\hline
\end{tabular}
\end{table}}

\footnotesize {\begin{table}[htbp]
\centering
\begin{tabular}{|c||c|c|c|c|c|c|c|   }
\hline
$case~n$ &2 & 3  & 4& 5     \\
\hline $~$ &&&& \\[-2.5ex]
$\pi^n(\Sigma^{n+7}\C P^2)$ & $2+3$ &$2^2+21$  &  $4+2^3+21$ & $4+2^2+63$   \\
\hline
$case~n$  & 6 &7&8 &9     \\
\hline  $~$ &&&& \\[-2.5ex]
$\pi^n(\Sigma^{n+7}\C P^2)$  &  $4+2^2+63$  & $8+2^2$  & $8+2^2$ &$8+2^2$   \\
\hline
$case~n$ &10 &11  & 12 &$n \geq 13$     \\
\hline  $~$ &&&& \\[-2.3ex]
$\pi^n(\Sigma^{n+7}\C P^2)$   & $\infty+8+2+63$ & $8+2+63$&  $\infty+8+2+63$   & $8+2+63$  \\
 \hline
\end{tabular}
\end{table}}

 {\begin{table}[htbp]
 \setlength\extrarowheight{2.5pt}
 \centering
\begin{tabular}{|c||c|c|c|c|c|c|c|c|c|c|c|}
\hline
$case~n$ &2 & 3  & 4& 5  & 6 &7 &8 \\
\hline
$~$ &&&&&&\\[-3.7ex]
$\pi^n(\Sigma^{n+8}\C P^2)$ & $2^2+21$ & $2+3 $& $8+2$ & $4+9 $ & $8+4+9+3 $ & $8+3 $ & $8^2+2^2$\\
\hline
$case~n$    &9 &10 &11  & 12 & 13 & $n \geq 14$&\\
\hline  $~$ &&&&&&\\[-3.6ex]
$\pi^n(\Sigma^{n+8}\C P^2)$   &$ 8+3$ &$ 8+2 $& $2^2 $& $ 2+3$ & 2+3 &3 &\\
\hline
\end{tabular}
\end{table}}

\normalsize{ }

\normalsize{ }

where an integer $n$ denotes the cyclic
group $\Z_n$,
$(s)^k$ denotes the $k$-times direct sum of $\Z_s$,
$\infty$ denotes the group of integers $\Z$ and ``+'' denotes the direct sum of abelian groups.
\\

In section 6, using our result in section 3, 4 and 5, we determine $n$th-homotopy groups of function spaces $map_\ast(\C P^2, \C P^2;\ast)$ for $4 \le n \le 13.$
\\

In Section 7, we apply our results to the classification of path components
of mapping spaces up to homotopy equivalent and evaluation fibrations up to fibre homotopy equivalent.
Hansen proved that the evaluation fibration
$w_f:map(\Sigma X,\Sigma Y;f) \to \Sigma Y$ has a section if and only
if $[f,id_{\Sigma Y}]=0$, where $[\,\,\,,\,\,]$ is the generalized
Whitehead product \cite{Hansen}.
Lupton and Smith proved that the
following statements are all equivalent:
(1) a map $f: X \to Y$ is cyclic,
(2) $w_f$ has a section, and
(3) two fibrations $w_f$ and
$w_\ast$ are fiber homotopy equivalent, where $\ast$ is a constant map \cite{LS}.
Moreover, we apply our results to the formulation of generalized Gottlieb groups
from suspended complex plane to sphere and Gottlieb groups of path components of constant map.
We use the notation of \cite{T,KMNST} freely.
\\

\section{Preliminaries}

 The complex projective plane $\C P^2$ is defined by the mapping cone
$S^2 \cup_{\eta_2} e^4$ where $\eta_2:S^3 \to S^2$ is the Hopf fibration.
We denote {$\eta_k=\Sigma^{k-2}\eta_2$ for $k\geq 2$}.
Consider a Puppe sequence
$$ S^3\xrightarrow{\eta_2}  S^2 \xrightarrow{i} \C P^2 \xrightarrow{p}
S^4 \xrightarrow{\eta_3} S^3 \xrightarrow{\Sigma i} \cdots,$$
where $i:S^2 \to \C P^2$ is the inclusion map and
$p:\C P^2 \to S^4$ is the collapsing map of $S^2$ to a point.
The Puppe sequence induces a long exact sequence of homotopy sets
$$\cdots \to \pi_{n+3}(S^m) \xrightarrow{\eta_{n+3}^\ast} \pi_{n+4}(S^m) \xrightarrow{\Sigma^n p^\ast}
[\Sigma^n \C P^2,S^m]$$
$$ \hspace{35mm}\xrightarrow{\Sigma^n i^\ast} \pi_{n+2}(S^m) \xrightarrow{\eta_{n+2}^\ast}
\hspace{2mm}\pi_{n+3}(S^m) \to \cdots.$$
Then we have the short exact sequence
\begin{equation}
\tag{1.1} 0 \to Coker \eta_{n+3}^\ast  \xrightarrow{\Sigma^n p^\ast}
[\Sigma^n \C P^2,S^m] \xrightarrow{\Sigma^n i^\ast} Ker \eta_{n+2}^\ast \to 0.
\end{equation}

Let $g_6(\C) : \Sigma^7 \C P^2 \to S^6$ be the $S^1$-transfer map \cite{M1993}. This map is the adjoint of the composite of inclusions
$$
\Sigma \C P^2 \hookrightarrow SU(3) \hookrightarrow SO(6) \hookrightarrow \Omega^6 S^6.
$$
We set
$
g_{n+6}(\C) =\Sigma^n g_6(\C)
$
for $n \ge 1$.
Note that
$
g_6(\C) = [\iota_6,\iota_6]\circ \Sigma^7 p +\nu_6 \circ \overline{2\iota_9}
$
and
$
2g_n(\C) = \nu_n \circ \overline{2\iota_{n+3}}
$
for $n \ge 7$ \cite[Proposition 3.3]{KMNST}.
\\

Let $G$ be an abelian group.
For a prime $p$, we
denote the $p$-primary parts of $G$ by $G_{(p)}$.
For an odd prime $p$, we have an isomorphism
\begin{equation}
\tag{1.2} [\Sigma^n \C P^2,S^k]_{(p)} \cong \pi_{n+2}(S^k)_{(p)} \oplus \pi_{n+4}(S^k)_{(p)},
\end{equation}
since $\pi_{n+1}(S^n)$ is of order 2 for $n \geq 3$
\cite[Proposition 5.1]{T}.\\

The space $\Sigma \C P^{n-1}$ has reduced homology
$$\widetilde{H}_\ast(\Sigma \C P^{n-1}) = \Z_p[x_3,x_5,\cdots,x_{2n-1}].$$
There is a Steenrod operation on $\widetilde{H}_\ast(\Sigma \C P^{n-1})$ given by
$$\mathcal{P}^j(x_{2r+1}) = r!/(j! \cdot(r-j)!) \cdot x_{2r+jp+1} \,\,\,\,(2 \le r \le n-1).$$
It has been shown that there exists a wedge decomposition \cite{MNT}
$$\Sigma \C P^{n-1} \simeq \bigvee_{i=1}^{p-1} A_i(n)$$
where
$A_i(n)$ has a reduced homology
$$
\widetilde{H}_\ast(A_i(n)) =\Z/p\Z\{x_{2i+1}, x_{2i+q+1}, \cdots, x_{2i+(r^n_i -1)q+1}\,\,|\,\,\, r^n_i  = \lfloor (n-i-1)/(p-1)\rfloor +1 \}
$$
that inherits a Steenrod operation from $\widetilde{H}_\ast(\Sigma \C P^{n-1}).$
\\


Consider elements $\alpha \in [Y,Z]$, $\beta \in [X,Y]$, and $\gamma \in [W,X]$ which satisfy $\alpha \circ \beta =0$ and $\beta \circ \gamma =0$.
Let $C_{\beta}$ be the mapping cone of $\beta$, $i:Y \to C_\beta$ and $p:C_\gamma \to \Sigma X$ be the inclusion and the shrinking map, respectively.
We denote  an {\it extension} $\overline{\alpha} \in [C_\beta ,Z]$ of $\alpha$ satisfying
$i^\ast(\overline{\alpha})=\alpha$
and  a {\it coextension} $\widetilde{\gamma} \in [\Sigma W ,C_\beta]$ of $\gamma$ satisfying
$p_\ast(\widetilde{\gamma})=\Sigma \gamma$ \cite{T}.
From the definition of extension, an extension $\overline{\alpha}$ exists if and only if $\alpha \circ \beta =0$. Similarly, a coextension $\widetilde{\gamma}$ exists if and only if $\beta \circ \gamma = 0.$
\\

We recall a useful relation between extensions and Toda brackets \cite[Proposition 1.9]{T}:
\begin{prop}
Let $\alpha \in [Y,Z]$, $\beta \in [X,Y]$, and $\gamma \in [W,X]$
be elements which satisfy $\alpha \circ \beta=0$ and $\beta \circ \gamma=0$.
Let $\{\alpha,\beta,\gamma\}$ be the Toda bracket and
$p:C_\gamma \to \Sigma W$ be the shrinking map.
Then, we have
$$
\alpha \circ \overline{\beta} \in \{\alpha,\beta,\gamma\} \circ p.
$$
\end{prop}

%
%

Here, we recall the concept of a cyclic map and Gottlieb group  $G_n(X)$ of
a space $X$ \cite{Got,Vara}.

\begin{definition}
A map $f :Y \to X$ is cyclic if there is a map $F:X \times Y \to
X$, called an affiliated map of $f$, such that the diagram homotopy commutative
$$\xymatrix{
X \times Y \ar[rr]^F && X  \\%
X \vee Y \ar[rr]^{1 \vee f} \ar[u]^{i} && X \vee X
\ar[u]_\nabla.}$$
\end{definition}

Let $G(Y,X)$ denote the set of all homotopy classes of cyclic maps from $Y$ to $X$.
Varadarajan showed that $G(Y,X)$ has group structure for any co-H-space $Y$ \cite{Vara}.
For an integer $n \ge 1$, the set of homotopy classes of cyclic maps $\Sigma^n X \to Y$ we denote by $G_n(X,Y)$, and call the {\it $n$-th generalized Gottlieb group of $(X,Y)$}.
When $Y=S^n$, $G(Y,X)=G_n(X)$ is the {\it $n$-th Gottlieb group of} $X$.
In \cite{Got}, Gottlieb introduced and studied the evaluation subgroups
$$G_n(X) =w_\ast(\pi_n(map(X,X;1)))$$
where $w : \ map(X,X;1)  \to X$ is the evaluation map.
Note that the $G_n(X)$ can alternatively be described as homotopy classes of maps
$f:S^n \to X$ such that
$(f|1):S^n \vee X \to X$ admits an extension $F:S^n \times X \to X$ up to homotopy.

\section{Determination of $[\Sigma^{n+6} \C P^2,S^n]$ }

In this section, we shall determine the generators of the 2-primary components of $n$-th cohomotopy groups of $(n+6)$-fold
suspended complex projective planes.
We analyze the following extension
\begin{equation}
\tag{2.1} 0 \to Coker\,\eta^\ast_{n+9} \xrightarrow{\Sigma^n p^\ast} [\Sigma^{n+6} \C P^2,S^n ] \xrightarrow{\Sigma^n i^\ast} Ker\, \eta^\ast_{n+8} \to 0
\end{equation}
induced by (1.1).
\begin{lem}
(1) For the homomorphisms $\eta^\ast_{n+9} : \pi^n_{n+9} \to \pi^n_{10}$, we have the following table of the cokernel of $\eta^\ast_{n+9}:$
{\begin{table}[htbp]
\setlength\extrarowheight{2.5pt}
\centering
\begin{tabular}{|c|c|c|c|c|c|c|c| c|c|c|c|c|c|}
\hline
$n=$ &2 & 3  & 4& 5, 6, 7 &8&9&10&11& $ \ge 12$\\
\hline
$~$ &&&&&\\[-3.7ex]
$Coker\,\eta^\ast_{n+9}  $ & $\Z_2$ & $Z_2$& $\Z_8 \oplus \Z_2$ & $\Z_4$   & $\Z_4^2$ & $Z_4$ & $\Z_2$ & $\Z_2$ & 0\\
\hline
generators  & $ \eta_2\mu_3 $ & $\e'$  & $\nu_4\s' $ & $\nu_n\s_{n+3}$  & $\s_8\nu_{15} $ & $\s_9\nu_{16}$ & $\s_{10}\nu_{17}$ & $\s_{11}\nu_{18}$   & 0     \\
   &   &   & $ \Sigma\e'$ &    & $ \nu_8\s_{11}$ &  &   &     & 0    \\
\hline
\end{tabular}
\end{table}}
\\

(2) The homomorphisms $\eta_{n+8}^n : \pi^n_{n+8}\to \pi^n_{n+9}$ are monomorphic except $n=6$ where
$$
Ker \,\eta^\ast_{14} = \Z_4\{2\bar{\nu}_6\}.
$$
\end{lem}

From (2.1) and Lemma 1, we state our result.
\begin{prop}
(1) $[\Sigma^{8} \C P^2,S^{2}]
= \Z_2\{\eta_2 \circ \mu_3 \circ \Sigma^8p\} \oplus \Z_{15}$.\\
(2) $[\Sigma^{9} \C P^2,S^{3}]
= \Z_2\{\e' \circ \Sigma^9 p\} \oplus \Z_3$.\\
(3) $[\Sigma^{10} \C P^2,S^{4}]
= \Z_8\{\nu_4 \circ \s' \circ \Sigma^{10}p\} \oplus \Z_2\{\Sigma\e'\circ \Sigma^{10}p\} \oplus \Z_3^2 \oplus \Z_5.$\\
(4) $[\Sigma^{11} \C P^2,S^{5}]
=  \Z_4\{ \nu_5 \circ \s_8\circ \Sigma^{11}p\} \oplus \Z_9$.\\
(5) $[\Sigma^{12} \C P^2,S^{6}]
=  \Z_4^2\{ \nu_6 \circ \s_9\circ \Sigma^{12} p, \bar{\nu}_6\circ \overline{2\iota_{14}}\} \oplus \Z_9 \oplus \Z_3$.\\
(6) $[\Sigma^{13} \C P^2,S^{7}]
=  \Z_4\{ \nu_7 \circ \s_{10}\circ \Sigma^{13} p\} \oplus \Z_3$.\\
(7) $[\Sigma^{14} \C P^2,S^{8}]
=   \Z_4^2\{  \s_8 \circ \nu_{15}\circ \Sigma^{14} p,
 \nu_8 \circ \s_{11}\circ \Sigma^{14} p\} \oplus \Z_3^2$.\\
(8) $[\Sigma^{15} \C P^2,S^{9}]
=   \Z_4\{ \s_9 \circ \nu_{16}\circ \Sigma^{15} p\} \oplus \Z_3$.\\
(9) $[\Sigma^{16} \C P^2,S^{10}]
=   \Z_2\{ \s_{10}\circ \nu_{17}\circ \Sigma^{16} p\} \oplus \Z_3$.\\
(10) $[\Sigma^{17} \C P^2,S^{11}]
=   \Z_2\{ \s_{11}\circ \nu_{18}\circ \Sigma^{17} p\} \oplus \Z_3$.\\
(11) $[\Sigma^{n+6} \C P^2,S^{n}]
=   \Z_3 $ for $n \geq 12$.
\begin{proof}
By Lemma 1 (2), we determine all group structures of $[\Sigma^{n+6}\C P^2,S^n]$ except $n=6$ so we consider $[\Sigma^{12}\C P^2,S^6]$.
By (2.1) and Lemma 1, we have the following commutative ladder:
$$\xymatrix{
0 \ar[r] & \Z_4\{\nu_6 \circ \s_9\} \ar[r]^(0.47){\Sigma^{12} p^\ast } \ar[d]^{\Sigma_1} &
[\Sigma^{12} \C P^2,S^{6}]_{(2)} \ar[r]^(0.6){\Sigma^{12}i^\ast } \ar[d]^{\Sigma_2} &
\Z_4\{2\overline{\nu}_6\} \ar[r] \ar[d]^{\Sigma_3} & 0 \\%
0 \ar[r] & \Z_4\{\nu_7 \circ \s_{10}\} \ar[r]^(0.48){\Sigma^{13}p^\ast } &
[\Sigma^{13} \C P^2,S^{7}]_{(2)} \ar[r]^(0.7){\Sigma^{12}i^\ast } &
0 \ar[r] & 0}$$
Since $\Sigma_1$ and $\Sigma^{13}p^\ast$ are isomorphisms,
$\Sigma^{12}p^\ast$ has left inverse.
This implies that the first row splits.
This completes the proof.
\end{proof}
\end{prop}

For odd prime $p$, we have an isomorphism
 \begin{equation}
\tag{2.2} [\Sigma^{n+6} \C P^2,S^n]_{(p)} \cong \pi_{n+8}(S^n)_{(p)} \oplus \pi_{n+10}(S^n)_{(p)}
 \end{equation}
by (2.2).
From \cite[Chapter 13]{T}, then, we have the following:

\begin{prop}
The odd primary components of $[\Sigma^{n+6} \C P^2,S^n]$ are\\
(1) $[\Sigma^8 \C P^2, S^2]  = \Z_3\{\eta_2 \circ \overline{\alpha_2(3)}\} \oplus \Z_5\{\eta_2 \circ \overline{\alpha_{1,5}(3)}\}.$\\
(2) $[\Sigma^9 \C P^2, S^3]  = \Z_3\{ \alpha_1(3)\circ\alpha_2(6) \circ \Sigma^9p \}$.\\
(3) $[\Sigma^{10} \C P^2, S^4]  = \Z_3^2\{ \alpha_1(4)\circ\alpha_2(7) \circ \Sigma^{10} p ,\nu_4 \circ \alpha_2(7) \circ \Sigma^{10}p \}$.\\
(4) $[\Sigma^{11} \C P^2, S^5]  = \Z_9\{ \beta_1(5) \circ \Sigma^{11} p\}$.\\
(5) $[\Sigma^{12} \C P^2, S^6]  = \Z_9\{ \beta_1(6) \circ \Sigma^{12} p\} \oplus \Z_3\{[\iota_6,\iota_6]\circ \overline{\alpha_1(11)}\}$.\\
(6) $[\Sigma^{13} \C P^2, S^7]  = \Z_3\{ \beta_1(7) \circ \Sigma^{13} p\}  $.\\
(7) $[\Sigma^{14} \C P^2, S^8]  = \Z_3^2\{ \beta_1(8) \circ \Sigma^{14} p, [\iota_8,\iota_8]\circ  {\alpha_1(15)} \circ \Sigma^{14} \}  $.\\
(8) $[\Sigma^{n+6} \C P^2, S^n]  = \Z_3\{ \beta_1(n) \circ \Sigma^{n+6} p\}$ for $n \ge 8$.
\end{prop}

\section{Determination of $[\Sigma^{n+7} \C P^2,S^n]$}

In this section, we shall determine the generators of the 2-primary components of $n$-th cohomotopy groups of $(n+7)$-fold
suspended complex projective planes.
We analyze the short exact sequence
\begin{equation}
\tag{3.1} 0 \to Coker\,\eta^\ast_{n+10} \xrightarrow{\Sigma^n p^\ast} [\Sigma^{n+7} \C P^2,S^n ] \xrightarrow{\Sigma^n i^\ast} Ker\, \eta^\ast_{n+9} \to 0
\end{equation}
induced by (1.1).
\begin{lem}
(1) For the homomorphisms $\eta^\ast_{n+10} : \pi^n_{n+10} \to \pi^n_{11}$, we have the following table of the cokernel of $\eta^\ast_{n+10}:$
{\begin{table}[htbp]
\setlength\extrarowheight{2.5pt}
\centering
\begin{tabular}{|c|c|c|c|c|c|c|c| c|c| c|c|c|c|}
\hline
$n=$ &2 & 3  & 4& 5  & 6, 7, 8, 9 & 10, 11 & 12 &$\ge 13$ \\
\hline
$~$ &&&&&\\[-3.7ex]
$Coker\,\eta^\ast_{n+10}  $ & $\Z_2$ & $Z_2^2$& $\Z_2^4$ & $\Z_4 \oplus \Z_2$  & $\Z_4 \oplus \Z_2$ & $\Z_4$ & $\Z \oplus \Z_4$& $\Z_4$  \\
\hline
generators  & $ \eta_2\e' $ & $\mu' $         & $\Sigma\mu',\nu_4\bar{\nu}_7 $ & $\zeta_5 $         & $\zeta_n $ & $ \zeta_n $   & $P(\iota_{25}) $  & $\zeta_n$  \\
            &               & $\e_3\nu_{11}$  & $\nu_4\e_7, \e_4\nu_{12}$ & $ \nu_5\bar{\nu}_8$  & $ \bar{\nu}_n\nu_{n+{8}}$  &  &$\zeta_{12}$ & \\
\hline
\end{tabular}
\end{table}}
\\

(2) For the homomorphisms $\eta^\ast_{n+9} : \pi^n_{n+9} \to \pi^n_{10}$, we have the following table of the kernel of $\eta^\ast_{n+9}:$
{\begin{table}[htbp]
\setlength\extrarowheight{2.5pt}
\centering
\begin{tabular}{|c|c|c|c|c|c|c|c|c|c|c|c|c|     }
\hline
$n=$ &2, 3  & 4, 5, 6 & 7  & 8 & 9 & 10 & $ \ge  11$  \\
\hline
$~$ &&&&\\[-3.7ex]
$Ker\,\eta^\ast_{n+9}  $ &   0 & $\Z_2$ & $\Z_2^2$  & $\Z_2^2$ & $\Z_2^2$ & $\Z \oplus \Z_2^2$ & $\Z_2^2$   \\
\hline
generators     &   0 & $\nu^3_n$ & $\nu^3_7 $                     & $\nu^3_8 $                       & $ \nu^3_9  $   & $2P(\iota_{21}) $  &  $\nu^3_n $    \\
               &     &           & $\s'\eta^2_{14}+\eta_7\e_8$   & $ (E\s')\eta^2_{15}+\eta_8\e_9$  & $ \eta_9\e_{10} $   & $ \nu^3_{11}, \eta_{11}\e_{12}$  &  $ \eta_n\e_{n+1}$   \\
\hline
\end{tabular}
\end{table}}

\end{lem}

By Lemma 2, we obtain
$$
\Sigma^n p^\ast:Coker\,\eta^\ast_{n+10} \to [\Sigma^{n+7}\C P^2,S^n]
$$
are isomorphic for $n=2, 3$.
So, we have the following:
\begin{prop}
(1) $[\Sigma^{9} \C P^2,S^{2}]
=  \Z_2\{ \eta_2 \circ \e'\circ \Sigma^9 p\} \oplus \Z_3$.\\
(2) $[\Sigma^{10} \C P^2,S^{3}]
= \Z_2^2\{  \mu'\circ \Sigma^{10} p,\e_3 \circ \nu_{11}\circ \Sigma^{10} p \} \oplus \Z_3 \oplus \Z_7 $.
\end{prop}

\begin{prop}
$[\Sigma^{11} \C P^2,S^{4}]
= \Z_4\{\nu_4^2 \circ g_{10}(\C)\} \oplus \Z_2^3\{(\Sigma\mu') \circ \Sigma^{11} p,\e_4 \circ \nu_{12} \circ \Sigma^{11} p,
\nu_4 \circ \overline{\nu}_7 \circ \Sigma^{11} p\} \oplus \Z_{21} $.\\
Relation : $\nu_4 \circ \e_7 \circ \Sigma^{11}p= 2\nu_4^2 \circ   g_{10}(\C).$
\begin{proof}
By (3.1) and Lemma 2, we have
$$0 \to  \Z_2^4\{  \Sigma\mu',\nu_4 \circ \overline{\nu}_7,\nu_4 \circ \e_7,\e_4 \circ \nu_{12}\} \oplus \Z_{15}
 \xrightarrow{\Sigma^{11} p^\ast} [\Sigma^{11} \C P^2,S^{4}]\xrightarrow{\Sigma^{11} i^\ast}\Z_2\{\nu_4^3\}\to 0 .$$
By \cite[Proposition 3.6]{KMNST}, 
we have a relation
$$2\nu_4^2 \circ g_{10}(\C) = \nu_4^3 \circ\overline{2\iota_{13}} = \nu_4 \circ\e_7\circ \Sigma^{11} p.$$
This completes the proof.
\end{proof}
\end{prop}

\begin{prop}
$[\Sigma^{12} \C P^2,S^{5}]
= \Z_4\{ \zeta_5\circ \Sigma^{12}p\} \oplus \Z_2^2\{\nu_5 \circ \overline{\nu}_8\circ \Sigma^{12}p,\nu_5^2\circ g_{11}(\C)\} \oplus \Z_7 \oplus \Z_9
$.
\begin{proof}
By (3.1) and Lemma 2, we have
$$
0 \to  \Z_4\{ \zeta_5\} \oplus \Z_2\{\nu_5 \circ \overline{\nu}_8\} \oplus \Z_7 \oplus \Z_9
 \xrightarrow{\Sigma^{12} p^\ast}
[\Sigma^{12} \C P^2,S^{5}]
\xrightarrow{\Sigma^{12} i^\ast}
\Z_2\{\nu_5^3\}
\to 0
$$
Since $\nu^2_5$ is of order 2,
$
\nu^2_5\circ g_{11}(\C) \in [\Sigma^{12}\C P^2,S^5]
$
is of order 2.
This completes the proof.
\end{proof}
\end{prop}

\begin{prop}
$[\Sigma^{13} \C P^2,S^{6}]
= \Z_4\{ \zeta_6\circ \Sigma^{13}p\} \oplus \Z_2^2\{ \overline{\nu}_6 \circ \nu_{14}\circ \Sigma^{13}p,\nu_6^2 \circ g_{12}(\C)\} \oplus \Z_{7} \oplus \Z_9 $.
\begin{proof}
By (3.1) and Lemma 2, we have
$$ 0 \to  \Z_4\{ \zeta_6\} \oplus \Z_2\{  \overline{\nu}_6 \circ \nu_{14}\} \oplus \Z_7 \oplus \Z_9
 \xrightarrow{\Sigma^{13} p^\ast}
[\Sigma^{13} \C P^2,S^{6}]
\xrightarrow{\Sigma^{13} i^\ast}
\Z_2\{\nu_6^3\}
\to 0 .$$
Then we have a commutative diagram
$$\xymatrix{
0 \ar[r] &  \Z_4\{ \zeta_5\} \oplus \Z_2\{\nu_5 \circ \overline{\nu}_8\}\ar[r]^(0.6){\Sigma^{12} p^\ast} \ar[d]^{\Sigma_1} & [\Sigma^{12} \C P^2,S^{5}] \ar[r]^(0.6){\Sigma^{12} i^\ast} \ar[d]_{\Sigma_2} &
\Z_2\{\nu_5^3\} \ar[r] \ar[d]^{\Sigma_3} &  0\\%
0 \ar[r] &  \Z_4\{ \zeta_6\} \oplus \Z_2\{ \overline{\nu}_6 \circ \nu_{14}\}\ar[r]^(0.61){\Sigma^{13} p^\ast} & [\Sigma^{13} \C P^2,S^{6}] \ar[r]^(0.59){\Sigma^{13} i^\ast}   &
\Z_2\{\nu_6^3\} \ar[r]  &  0
}$$
%
Since the first row is split and $\Sigma_3$ is an isomorphism, the second row is also split.
This completes the proof.
\end{proof}
\end{prop}

Note that the Hopf map $\s_8:S^{15} \to S^8$ induces an isomorphism

\begin{equation}[X,S^7] \oplus [\Sigma X,S^{15}] \to [\Sigma X,S^8], \,\,\,\,(\alpha,\beta) \mapsto \Sigma \alpha + \s_8 \circ \beta.
\end{equation}

\begin{lem}
$$[\Sigma^{14}\C P^2,S^7]_{(2)} \cong [\Sigma^{15}\C P^2,S^8]_{(2)} \cong [\Sigma^{16}\C P^2,S^9]_{(2)}.$$
\begin{proof}
By (2) and \cite[Proposition 3.2]{KMNST}, we have an isomorphism
$$\Sigma :[\Sigma^{14}\C P^2,S^7] \to[\Sigma^{15}\C P^2,S^8] .$$
We consider a 2-primary EHP sequence
$$
[\Sigma^{17}\C P^2,S^{17}]_{(2)} \xrightarrow{\triangle } [\Sigma^{15}\C P^2,S^{8}]_{(2)} \,\,\xrightarrow{\Sigma  }
[\Sigma^{16}\C P^2,S^{9}]_{(2)} \xrightarrow{ H}$$
$$\hspace{21mm} \xrightarrow{ H}[\Sigma^{16}\C P^2,S^{17}]_{(2)} \xrightarrow{\triangle }
[\Sigma^{14}\C P^2,S^{8}]_{(2)}
$$
where $[\Sigma^{17}\C P^2,S^{17}] =0$ ,
$[\Sigma^{16}\C P^2,S^{17}]_{(2)}=\Z_4\{\nu_{17} \circ \Sigma^{16}p\}$
and
$[\Sigma^{14}\C P^2,S^{8}]_{(2)} = \Z_4^2\{\s_8 \circ \nu_{15} \circ \Sigma^{14}p,\nu_8 \circ \s_{11} \circ \Sigma^{14}p\} $
by Proposition 3 and \cite[Proposition 3.1, 3.2]{KMNST}.
Then $\Sigma: [\Sigma^{15}\C P^2,S^{8}]_{(2)} \to  [\Sigma^{16}\C P^2,S^{9}]_{(2)}$ is injective.
By \cite[(7.19)]{T}, we have
$$\Delta(\nu_{17} \circ \Sigma^{16}p)
= P(\nu_{17}) \circ \Sigma^{14}p
= 2\s_8  \circ \nu_{15}\circ \Sigma^{14}p -x \nu_8 \circ \s_{11} \circ \Sigma^{14}p$$
for an odd $x$.
Then $\triangle : [\Sigma^{16}\C P^2,S^{17}]_{(2)} \to [\Sigma^{14}\C P^2,S^{8}]_{(2)}$
is injective, so that
$\Sigma: [\Sigma^{15}\C P^2,S^{8}]_{(2)} \to  [\Sigma^{16}\C P^2,S^{9}]_{(2)}$
is surjective.
Therefore $\Sigma$ is isomorphic.
This completes the proof.
\end{proof}
\end{lem}

By \cite[(5.5), (7.10), (7.19)]{T}, we have
$$
\s' \circ \eta^3_{14} = \s' \circ 4\nu_{14} =4(\nu_7 \circ \s_{10})= \eta^2_7 \circ \e_9 =\eta_7 \circ \e_8 \circ \eta_{16}.
$$
This implies $\s' \circ \eta^3_{14} + \eta_7  \circ \e_8 \circ \eta_{16}=0$.
Then we can take an {\it extension} $\overline{\s' \circ \eta^2_{14} + \eta_7  \circ \e_8} \in [\Sigma^{14} \C P^2,S^7]$ of $\s' \circ \eta^2_{14} + \eta_7  \circ \e_8$.
Similarly, we can take an {\it extension} $\overline{(E\s') \circ \eta_{15}^2 +\eta_8 \circ \e_9} \in [\Sigma^{15} \C P^2,S^{8}]$ of $ {(E\s') \circ \eta_{15}^2 +\eta_8 \circ \e_9}$.
By \cite[(7.5), (7.10), (7.19), Theorem 7.2]{T}, we have
$$
\eta_{n+9}^\ast(\eta_n\e_{n+1} ) = \eta^2_n\e_{n+3}=4\nu_n\s_{n+3} = 0
$$
for $n \ge9$.
Then we can take {\it extensions} $\overline{\eta_n \circ \e_{n+1}} \in [\Sigma^{n+7} \C P^2,S^{n}]$ of ${\eta_n \circ \e_{n+1}}$ for $n \ge 9$.

\begin{prop}
(1)
$[\Sigma^{14} \C P^2,S^{7}]
= \Z_8\{\overline{\s' \circ \eta_{14}^2 +\eta_7 \circ \e_8}\} \oplus \Z_2^2\{\overline{\nu}_7 \circ \nu_{15} \circ \Sigma^{14}p, \nu_7^2 \circ g_{13}(\C) \}$.\\
Relation: $2\overline{\s' \circ \eta_{14}^2 +\eta_7 \circ \e_8} = x \zeta_7 \circ \Sigma^{14}p$ for an odd $x$.\\
(2) $[\Sigma^{15} \C P^2,S^{8}]
=
\Z_8\{\overline{(\Sigma\s') \circ \eta_{15}^2 +\eta_8 \circ \e_9}\} \oplus \Z_2^2\{\overline{\nu}_8 \circ \nu_{16} \circ \Sigma^{15}p,\nu_8^2 \circ g_{14}(\C)\}$.\\
Relation: $2\overline{(\Sigma\s') \circ \eta_{15}^2 +\eta_8 \circ \e_9} = x \zeta_8 \circ \Sigma^{15}p$ for an odd $x$.\\
(3)
$[\Sigma^{16} \C P^2,S^{9}]
= \Z_8\{\overline{\eta_9 \circ \e_{10}}\} \oplus \Z_2^2\{\overline{\nu}_9 \circ \nu_{17} \circ \Sigma^{16}p,\nu_9^2 \circ g_{15}(\C)\}$.\\
Relation: $2 \overline{\eta_9 \circ \e_{10}} = x\zeta_9 \circ \Sigma^{16}p   $ for an odd $x$.
\begin{proof}
By (3.1) and Lemma 2, we have
$$0 \to \Z_4\{ \zeta_7\} \oplus \Z_2\{\overline{\nu}_7 \circ \nu_{15}\} \oplus \Z_7 \oplus \Z_9
 \xrightarrow{\Sigma^{14} p^\ast}
[\Sigma^{14} \C P^2,S^{7}]
\xrightarrow{\Sigma^{14} i^\ast}
\Z_2^2\{\nu_7^3,\s' \circ \eta_{14}^2+\eta_7 \circ \e_8\}
\to 0 $$
$$0 \to  \Z_4\{ \zeta_8\} \oplus \Z_2\{\overline{\nu}_8 \circ \nu_{16}\} \oplus \Z_7 \oplus \Z_9
 \xrightarrow{\Sigma^{15} p^\ast}
[\Sigma^{15} \C P^2,S^{8}]
\xrightarrow{\Sigma^{15} i^\ast}
\Z_2^2\{\nu_8^3,(\Sigma\s') \circ \eta_{15}^2+\eta_8 \circ \e_9\}
\to 0 $$
$$0 \to \Z_4\{ \zeta_9\} \oplus \Z_2\{\overline{\nu}_9 \circ \nu_{17}\} \oplus \Z_7 \oplus \Z_9
 \xrightarrow{\Sigma^{16} p^\ast}
[\Sigma^{16} \C P^2,S^{9}]
\xrightarrow{\Sigma^{16} i^\ast}
\Z_2^2\{\nu_9^3,\eta_9 \circ \e_{10}\}
\to 0 .$$
By \cite[Lemma 5.14]{T}, we have
$
\Sigma((\Sigma\s') \circ \eta_{15}^2) = (\Sigma^2 \s') \circ \eta_{16}^2 = (2\s_9) \circ \eta_{16}^2 =  0.
$
Then we obtain
$$
\Sigma( \overline{(\Sigma\s')\circ \eta_{15}^2+\eta_8 \circ \e_9} ) \equiv \overline{\eta_9 \circ \e_{10}} \,\,\bmod \,\, \zeta_9\circ\Sigma^{16}p, \bar{\nu}_9\circ \nu_{17}\circ \Sigma^{16}p.
$$
This implies that
$2 \overline{\eta_9 \circ \e_{10}} = 2 \iota_9 \circ \overline{\eta_9 \circ \e_{10}}$ by \cite[Theorem 6.7]{Hilton}.
By \cite[Proposition 1.9, Lemma 9.1]{T}, we obtain
$$2\iota_9 \circ \overline{\eta_9 \circ \e_{10}} \in \{2\iota_9,\eta_9\circ \e_{10},\eta_{18}\} \circ \Sigma^{16}p \ni \zeta_9 \circ \Sigma^{16}p
\,\,\, \bmod \,\, 2\zeta_9 \circ \Sigma^{16}p,\,\overline{\nu}_9 \circ \nu_{17} \circ \Sigma^{16} p$$
that is,
$$2 \overline{\eta_9 \circ \e_{10}} \equiv x\zeta_9 \circ \Sigma^{16}p \,\,\bmod \, \overline{\nu}_9 \circ \nu_{17} \circ \Sigma^{16} p$$
for an odd $x$.
Since
$2(\nu_9^2 \circ g_{15}(\C)) =2(\nu_9^2 \circ \Sigma g_{14}(\C))= (2\nu_9^2) \circ g_{15}(\C)=0$ by \cite[(2.1)]{T},
$\nu_9^2 \circ g_{15}(\C)$ is of order 2.
By Lemma 3, we complete the proof.

\end{proof}
\end{prop}

By \cite[(7.21)]{T}, we have
$
P(\iota_{21})\circ \eta_{19} = P(\eta_{21}) = 2\sigma_{10}\circ \nu_{17}.
$
Then we have
$
2P(\iota_{21})\circ \eta_{19} = 4 \s_{10}\circ \nu_{17}=0
$ \cite[Theorem 7.3]{T}.
Thus we can take an {\it extension} $P(\iota_{21}) \circ \overline{2\iota_{10}} \in [\Sigma^{17} \C P^2,S^{10}]$ of  $ {2P(\iota_{21})}$.

\begin{prop}
(1) $[\Sigma^{17} \C P^2,S^{10}]
= \Z\{P(\iota_{21}) \circ \overline{2\iota_{19}}\} \oplus \Z_8\{\overline{\eta_{10} \circ \e_{11}}\} \oplus \Z_2\{\nu_{10}^2 \circ g_{16}(\C)\} \oplus \Z_{7} \oplus \Z_9$.\\
(2) $[\Sigma^{18} \C P^2,S^{11}]
= \Z_8\{\overline{\eta_{11}\circ \e_{12}}\} \oplus \Z_2\{\nu_{11}^2 \circ g_{17}(\C)\} \oplus \Z_{7} \oplus \Z_9$.\\
(3) $[\Sigma^{19} \C P^2,S^{12}]
= \Z\{P(\iota_{25}) \circ \Sigma^{19}p\} \oplus \Z_8\{\overline{\eta_{12}\circ\e_{13}}\} \oplus \Z_2\{\nu_{12}^2 \circ g_{18}(\C)\} \oplus \Z_{7} \oplus \Z_9$.\\
(4) For $n \geq 13$, $[\Sigma^{n+7} \C P^2, S^n]
= \Z_8\{\overline{\eta_n\circ \e_{n+1}}\} \oplus \Z_2\{\nu_{n}^2 \circ g_{n+6}(\C)\}  \oplus \Z_{7} \oplus \Z_9 $.
\begin{proof}
(1) 
By (3.1) and Lemma 2, we have
$$0 \to \Z_4\{ \zeta_{10}\}  \oplus \Z_{7} \oplus \Z_9
 \xrightarrow{\Sigma^{17} p^\ast}
[\Sigma^{17} \C P^2,S^{10}]
\xrightarrow{\Sigma^{17} i^\ast}
\Z\{2P(\iota_{21})\} \oplus \Z_2^2\{\nu_{10}^3,\eta_{10} \circ \e_{11}\}
\to 0 .$$
By proposition 8, we have
$$
2 \overline{\eta_{10} \circ \e_{11}} \equiv x\zeta_{10} \circ \Sigma^{17}p \,\,\,\bmod\,\,\,\bar\nu_{10}\circ \nu_{18}\circ \Sigma^{17}p
$$
for an odd $x$.
Since $\nu_{11}^2$ is of order 2, we have $\nu_{11}^2 \circ g_{17}(\C)$ is of order 2.
Thus we have
$$
[\Sigma^{17} \C P^2,S^{10}]
= \Z\{\overline{2P(\iota_{21})}\} \oplus \Z_8\{\overline{\eta_{10} \circ \e_{11}}\} \oplus \Z_2\{\nu_{10}^2 \circ g_{16}(\C)\} \oplus \Z_{7} \oplus \Z_9 .
$$
(2)
By (3.1) and Lemma 2, we have
$$0 \to \Z_4\{ \zeta_{11}\}  \oplus \Z_7 \oplus \Z_9
 \xrightarrow{\Sigma^{18} p^\ast}
[\Sigma^{18} \C P^2,S^{11}]
\xrightarrow{\Sigma^{18} i^\ast}
 \Z_2^2\{\nu_{11}^3,\eta_{11} \circ \e_{12}\}
\to 0 .$$
By (1),
we obtain
$$
[\Sigma^{18} \C P^2,S^{11}]
= \Z_8\{\overline{\eta_{11}\circ \e_{12}}\} \oplus \Z_2\{\nu_{11}^2 \circ g_{17}(\C)\} \oplus \Z_{7} \oplus \Z_9 .
$$
Similarly, we can determine group structures of $[\Sigma^{n+7} \C P^2, S^n]$ for $n \ge 12.$
This completes the proof.
\end{proof}
\end{prop}

For odd prime $p$, we have an isomorphism
 \begin{equation}
[\Sigma^{n+7} \C P^2,S^n]_{(p)} \cong \pi_{n+9}(S^n)_{(p)} \oplus \pi_{n+11}(S^n)_{(p)}
 \end{equation}
by (2.2).
From \cite[Chapter 13]{T}, then, we have the following:

\begin{prop}
The odd primary components of $[\Sigma^{n+7} \C P^2,S^n]$ are\\
(1) $[\Sigma^9 \C P^2, S^2]  = \Z_3\{\eta_2 \circ  \alpha_3(3) \circ \Sigma^9 p \}  $\\
(2) $[\Sigma^{10} \C P^2, S^3]  = \Z_3\{ \alpha_3(3)\circ \Sigma^{10} p \} \oplus \Z_7\{\alpha_{1,7}(3) \circ \Sigma^{10}p\}$.\\
(3) $[\Sigma^{11} \C P^2, S^4]  = \Z_3\{ \alpha_3(4)\circ \Sigma^{11} p \} \oplus \Z_7\{\alpha_{1,7}(4) \circ \Sigma^{11}p\}$.\\
(4) $[\Sigma^{n+7} \C P^2, S^n]  = \Z_7\{\alpha_{1,7}(n) \circ \Sigma^{n+7}p\} \oplus \Z_9\{\alpha_3'(n)\}$ for $n \ge 5$.
\end{prop}

\section{Determination of $[\Sigma^{n+8} \C P^2,S^n]$}

In this section, we determine the generators of the 2-primary components of $n$-th cohomotopy groups of $(n+8)$-fold
suspended complex projective planes.
We analyze the short exact sequence
\begin{equation}
\tag{3.1} 0 \to Coker\,\eta^\ast_{n+11} \xrightarrow{\Sigma^n p^\ast} [\Sigma^{n+8} \C P^2,S^n ] \xrightarrow{\Sigma^n i^\ast} Ker\, \eta^\ast_{n+10} \to 0
\end{equation}
induced by (1.1).
\begin{lem}
(1) For the homomorphisms
$
\eta^\ast_{n+11} : \pi^n_{n+11} \to \pi^n_{n+12}
$,
 we have the following table of the cokernel of $\eta^\ast_{n+11}:$
{\begin{table}[htbp]
\setlength\extrarowheight{2.5pt}
\centering
\begin{tabular}{|c|c|c|c|c|c|c|c|c|c|c|c|c|c|c|   }
\hline
$n=$ & 2 & 3  & 4& 5  & 6 & 7, 8, 9 &  10 & 11&12 & 13 & $14  \ge $ \\
\hline
$~$ &&&&&\\[-3.7ex]
$Coker\,\eta^\ast_{n+11}  $ & $\Z_2^2$ & 0 & $\Z_2$ & 0  & $\Z_8$ &0 & $\Z_4$ & $\Z_2$& $\Z_2$& $\Z_2$ & 0  \\
\hline
generators  & $\eta_2\mu' $ & 0 & $\nu_4\mu_7$ & 0 & $P(\s_{13})$ & 0  & $P(\nu_{21})$ & $\theta'$ & $\theta$ & $\Sigma\theta$ & 0\\
           & $ \eta_2\nu'\bar{\nu}_6$ &   &   &  &   &    &   &  &   &   &  \\
\hline
\end{tabular}
\end{table}}
\\

(2) For the homomorphisms $\eta^\ast_{n+10} : \pi^n_{n+10} \to \pi^n_{n+11}$, we have the following table of the kernel of $\eta^\ast_{n+10}:$
{\begin{table}[htbp]
\setlength\extrarowheight{2.5pt}
\centering
\begin{tabular}{|c|c|c|c|c|c|c|c|c|c|c|c|c|c|      }
\hline
$n=$ &2 & 3  & 4 & 5, 6  & 7  & 8 & 9 & 10 & 11 & $12 \ge  $    \\
\hline
$~$ &&&&\\[-3.7ex]
$Ker\,\eta^\ast_{n+10}  $ & 0 & $\Z_2$ & $ \Z_4 \oplus \Z_2$ & $\Z_4$  & $\Z_8$ & $\Z_8^2$ & $\Z_8$ & $\Z_4$ & $\Z_2 $  &0   \\
\hline
generators      &0 & $\eta^2_3\e_5$ & $2(\nu_4\s') $ & $2\nu_n\s_{n+3}$  & $\nu_7\s_{10}$  & $\s_8\nu_{15} $ & $\s_9\nu_{16}$ & $\s_{10}\nu_{17}$ & $\s_{11}\nu_{18}$ & 0 \\
        &  &   & $ \eta^2_4\e_6$ &    &    & $ \nu_8\s_{11}$ &  &   &   &   \\
\hline
\end{tabular}
\end{table}}
\end{lem}
By Lemma 4, we have the following
\begin{prop}
$
[\Sigma^{10}\C P^2,S^2] = \Z_2^2\{\eta_2\circ \mu' \circ \Sigma^{10}p, \eta_2\circ \nu' \circ \bar{\nu}_6 \circ \Sigma^{10}p\} \oplus \Z_{21}.
$
\end{prop}
By \cite[(7.12)]{T}, we have
$
\e' \circ \eta_{13} = \nu' \circ \e_6
$
and $\nu' \circ \e_6
$
is of order 2.
Then we can take an extension
$
\overline{\eta^2_3\circ\e_5} \in [\Sigma^{11}\C P^2,S^3]
$
of $2\e'=\eta^2_3\circ \e_5$ \cite[Lemma 6.6]{T}.
Similarly, we can take an extension
$
\overline{\eta^2_4\circ\e_6} \in [\Sigma^{12}\C P^2,S^4]
$
of $ \eta^2_4\circ \e_6$.
Then we obtain the following.
\begin{prop}
$
[\Sigma^{11}\C P^2,S^3] = \Z_2\{\overline{\eta^2_3 \circ \e_5}\} \oplus \Z_3.
$
\end{prop}

\begin{prop}
$
[\Sigma^{12}\C P^2,S^4] = \Z_8\{\nu_4 \circ \s' \circ \overline{2\iota_{14}}\} \oplus \Z_2\{\overline{\eta^2_4\circ \e_6}\}.
$
\begin{proof}
By Lemma 4, we have a short exact sequence
$$
0 \to \Z_2\{\nu_4\circ \mu_7\} \xrightarrow{\Sigma^{12}p^\ast} [\Sigma^{12}\C P^2,S^4] \xrightarrow{\Sigma^{12}i^\ast} \Z_4\{2(\nu_4\circ\s')\} \oplus \Z_2\{2\Sigma\e'\} \to 0.
$$
By \cite[Proposition 3.7]{KMNST}, we obtain
$$
4\nu_4\circ \s' \circ \overline{2\iota_{14}} = \nu_4\circ \mu_7 \circ \Sigma^{12}p.
$$
This complete the proof.
\end{proof}
\end{prop}

By \cite[p.152]{T},
$
\nu_5\circ \s_8\circ \eta_{15}=\nu_5\circ \e_8
$
is of order 2.
Then we can take an extension
$
\nu_5\circ \s_8\circ \overline{2\iota_{15}} \in [\Sigma^{13}\C P^2,S^5]
$
of
$
2\nu_5\circ \s_8.
$
Similarly, we can take an extension
$
\nu_6\circ \s_9\circ \overline{2\iota_{16}} \in [\Sigma^{14}\C P^2,S^6]
$
of
$
2\nu_6\circ \s_9.
$
By Lemma 4, we have the following
\begin{prop}
(1)
$
[\Sigma^{13}\C P^2,S^5] = \Z_4\{\nu_5\circ \s_8 \circ \overline{2\iota_{15}}\} \oplus \Z_9.
$\\
(2)
$
[\Sigma^{14}\C P^2,S^6] = \Z_8\{P(\s_{13}) \circ \Sigma^{14}p\}
 \oplus \Z_4\{\nu_5\circ \s_8 \circ \overline{2\iota_{15}}\} \oplus \Z_9 \oplus \Z_3.
$
\begin{proof}
By Lemma 4, we have (1). So we only prove (2).
Consider the following commutative diagram induced by Lemma 4.
$$
\xymatrix{
0 \ar[r] & 0 \ar[r] \ar[d] & [\Sigma^{13}\C P^2,S^5] \ar[r]^{\Sigma^{13}i^\ast} \ar[d]^{\Sigma_1} & \Z_4\{2\nu_5\circ\s_8\} \ar[r] \ar[d]^{\Sigma_2}& 0 \\%
0 \ar[r] & \Z_8\{P(\s_{13})\} \ar[r]^{\Sigma^{13}p^\ast} & [\Sigma^{14}\C P^2,S^6] \ar[r]^{\Sigma^{14}i^\ast} & \Z_4\{2\nu_6\circ\s_9\} \ar[r] & 0
}$$
Since $\Sigma^{13}i^\ast$ and $\Sigma_2$ are isomorphic, the second row is split.
This complete the proof.
\end{proof}
\end{prop}

By \cite[(5.9)]{T},
$
\nu_n\circ \s_{n+3}\circ \eta_{n+10} =0
$
for
$n \ge 7$.
Then we can take an extension
$
\overline{\nu_n\circ \s_{n+3}} \in [\Sigma^{n+8}\C P^2,S^n]
$
of
$
\nu_n\circ \s_{n+3}.
$
Similarly, we have
$
\s_m\circ \nu_{m+7} \circ \eta_{m+10} = \nu_m\circ \s_{m+3}\circ \eta_{m+10}=0
$
for $m \ge 8.$
Then we can take extensions
$
\overline{\s_m\circ \nu_{m+7}}, \overline{\nu_m\circ \s_{m+3}} \in [\Sigma^{m+8} \C P^2,S^m]
$
of
$
\s_m\circ \nu_{m+7}
$
and
$
\nu_m \circ \s_{m+3}
$
respectively.
By Lemma 4, we have the following.
\begin{prop}
(1)
$
[\Sigma^{15}\C P^2,S^7] = \Z_8\{\overline{\nu_7\circ\s_{10}}\} \oplus \Z_3.
$\\
(2)
$
[\Sigma^{16}\C P^2,S^8] = \Z_8^2\{\overline{\s_8\circ \nu_{15}}, \overline{\nu_8\circ \s_{11}}\} \oplus \Z_3^2.
$\\
(3)
$
[\Sigma^{17}\C P^2,S^9] = \Z_8\{\overline{\s_9\circ \nu_{16}}\} \oplus \Z_3.
$
\end{prop}

\begin{lem}
$
\{2\s_{10}\circ \nu_{17},2\iota_{20},\eta_{20}\}_1 = 2P(\nu_{21}).
$
\begin{proof}
We recall that
$
\pi^{10}_{22}=\Z_4\{P(\nu_{21})\}
$
\cite[Proposition 7.6]{T}.
By \cite[Proposition 2.6, (7.21)]{T}, we have
$$
H\{2\s_{10}\circ \nu_{17}, 2\iota_{20},\eta_{20}\}_1 = \eta^3_{19}=4\nu_{19}.
$$
Since
$
(H \circ P)(\nu_{21}) =\pm2\nu_{19}
$
\cite[Proposition 2.7]{T},
we obtain the result.

\end{proof}
\end{lem}

\begin{prop}
$
[\Sigma^{18}\C P^2,S^{10}] = \Z_8\{\s_{10}\circ g_{17}(\C)\} \oplus
\Z_2\{2\s_{10}\circ g_{17}(\C) -P(\nu_{21})\circ \Sigma^{18}p\}.
$
\begin{proof}
By Lemma 4, we have a short exact sequence
$$
0 \to \Z_4\{P(\nu_{21})\} \xrightarrow{\Sigma^{18}p^\ast} [\Sigma^{18}\C P^2,S^{10}]
 \xrightarrow{\Sigma^{18}i^\ast} \Z_4\{\s_{10}\circ \nu_{17}\} \to 0.
 $$
By \cite[Proposition 1.9]{T}, we have
\begin{eqnarray}
\notag 4\s_{10}\circ g_{17}(\C) &=& 2\s_{10}\circ \nu_{17}\circ \overline{2\iota_{20}} \\
\notag &\in& \{2\s_{10}\circ \nu_{17}, 2\iota_{20}, \eta_{20}\}_1\circ \Sigma^{18}p\\
\notag &=& 2P(\nu_{21})\circ \Sigma^{18}p.
\end{eqnarray}
This completes the proof.
\end{proof}
\end{prop}

\begin{prop}
$
[\Sigma^{19}\C P^2,S^{11}] = \Z_2^2\{\s_{11}\circ g_{18}(\C), \theta'\circ \Sigma^{19}p\}.
$
\begin{proof}
By Lemma 4, we have a short exact sequence
$$
0 \to \Z_2\{\theta'\} \xrightarrow{\Sigma^{19}p^\ast} [\Sigma^{19}\C P^2,S^{11}]
 \xrightarrow{\Sigma^{19}i^\ast} \Z_2\{\s_{11}\circ \nu_{18}\} \to 0.
 $$
We consider 2-primary EHP-sequence
$$
[\Sigma^{20}\C P^2,S^{21}] \xrightarrow{P} [\Sigma^{18}\C P^2,S^{10}] \xrightarrow{\Sigma } [\Sigma^{19}\C P^2,S^{11}] \xrightarrow{H}
$$
$$
\hspace{16.5mm} \xrightarrow{H} [\Sigma^{19}\C P^2,S^{21}] \xrightarrow{P } [\Sigma^{17}\C P^2,S^{10}]
$$
where
$
[\Sigma^{20}\C P^2,S^{21}] = \Z_4\{\nu_{21}\circ \Sigma^{20}p\}
$
and
$
[\Sigma^{19}\C P^2,S^{21}] = \Z\{\overline{2\iota_{21}}\}
$
\cite[p.110, Proposition 3.1]{KMNST}.
Then
$
P(\nu_{21}\circ \Sigma^{20}p)=P(\nu_{21})\circ \Sigma^{18}p
$
and
$
P(\overline{2\iota_{21}}) \equiv P(\iota_{21})\circ \overline{2\iota_{19}}\,\,\,\bmod\,\,\,\nu^2_{10}\circ g_{16}(\C), 2\overline{\eta_{10}\circ \e_{11}}.
$
So, two P-homomorphisms are monomorphisms.
Thus
$
[\Sigma^{19}\C P^2,S^{11}] \cong \Z_2^2.
$
This completes the proof.
\end{proof}
\end{prop}
By Lemma 4, we have the following.
\begin{prop}
(1)
$
[\Sigma^{20}\C P^2,S^{12} ] = \Z_2\{\theta \circ \Sigma^{20}p\} \oplus \Z_3.
$\\
(2)
$
[\Sigma^{21}\C P^2,S^{13} ] = \Z_2\{(\Sigma\theta) \circ \Sigma^{21}p\} \oplus \Z_3.
$\\
(3)
$
[\Sigma^{n+8}\C P^2,S^{n} ] = \Z_3
$
for $n \ge 14.$
\end{prop}

For odd prime $p$, we have an isomorphism
 \begin{equation}
[\Sigma^{n+8} \C P^2,S^n]_{(p)} \cong \pi_{n+10}(S^n)_{(p)} \oplus \pi_{n+12}(S^n)_{(p)}
 \end{equation}
by (2.2).
From \cite[Chapter 13]{T}, then, we have the following:

\begin{prop}
The odd primary components of $[\Sigma^{n+8} \C P^2,S^n]$ are\\
(1) $[\Sigma^{10} \C P^2, S^2]  = \Z_3\{\eta_2 \circ  {\alpha_3(3)} \circ \Sigma^{10}p \} \oplus \Z_5\{\eta_2\circ {\alpha_{1,7}(3)}\circ \Sigma^{10}p\}  $.\\
(2) $[\Sigma^{11} \C P^2, S^3]  = \Z_3\{ \alpha_1(3)\circ \overline{\alpha_2(6)} \}$.\\
(3) $[\Sigma^{12} \C P^2, S^4]  = \Z_3^2\{ \alpha_1(3)\circ \overline{\alpha_2(6)} , [\iota_4,\iota_4]\circ \overline{\alpha_2(7)} \}$.\\
(4) $[\Sigma^{13} \C P^2, S^5]  = \Z_9\{\overline{\beta_1(5)}  \}$.\\
(5) $[\Sigma^{14} \C P^2, S^6]  = \Z_9\{\overline{\beta_1(6) } \} \oplus \Z_3\{[\iota_6,\iota_6]\circ \alpha_2(11)\circ \Sigma^{14}p\} \oplus \Z_5\{[\iota_6,\iota_6]\circ \alpha_{1,5}(11)\circ \Sigma^{14}p\}.$\\
(6) $[\Sigma^{15} \C P^2, S^7]  = \Z_3\{\overline{\beta_1(7)} \}.$\\
(7) $[\Sigma^{16} \C P^2, S^8]  = \Z_3^2\{\overline{\beta_1(8) }, [\iota_8,\iota_8]\circ \overline{\alpha_1(15)}\}.$\\
(8) $[\Sigma^{17} \C P^2, S^9]  = \Z_3\{\overline{\beta_1(9)}\}$.\\
(9) $[\Sigma^{18} \C P^2, S^{10}]  = \Z_3^2\{\overline{\beta_1(9)} , [\iota_{10},\iota_{10}]\circ \alpha_1(19)\circ \Sigma^{18}p\}$.\\
(10) $[\Sigma^{n+8} \C P^2, S^{n}]  = \Z_3\{\overline{\beta_1(n)}\}$ for $n \ge 11$.

\end{prop}


\section{Homotopy groups of $Map_\ast(\C P^2,\C P^2,\ast)$}

The complex projective plane $\C P^2$ is a base space of a $S^1$-bundle
$
S^1 \to S^5 \xrightarrow{p} \C P^2.
$
So, we have an isomorphism
$$p_\ast : [\Sigma^n \C P^2, S^5 ]\to [\Sigma^n \C P^2,\C P^2].$$
From the adjointness \cite{Lang}, we obtain an isomorphism
$$\pi_n(map_\ast(\C P^2,\C P^2;\ast)) \cong [\Sigma^n \C P^2,\C P^2].$$
McGibbon \cite{Mc} showed that homotopy class $[\C P^n, \C P^n] \cong \Z$
which is determined by using homomorphism between homology groups.
We denote the $n$-th homotopy group $\pi_n(Map_\ast(\C P^2,\C P^2;\ast))$ by $\pi_n(map_\ast(\C P^2,\C P^2 )).$
By the results of section 3, 4, 5 and \cite{KMNST}, we obtain the following:

\begin{theorem}
$
\pi_4(map_\ast(\C P^2,\C P^2)) \cong \Z_4\{\nu_5 \circ \Sigma^4p\} \oplus \Z_3\{\alpha_1(5) \circ \Sigma^4p\}
$. \\
$
\pi_5(map_\ast(\C P^2,\C P^2))  = 0
$. \\
$
\pi_6(map_\ast(\C P^2,\C P^2)) \cong \Z_4\{\nu_5 \circ \overline{2\iota_8}\} \oplus \Z_3\{\overline{\alpha_1(5)}\}.
$ \\
$
\pi_7(map_\ast(\C P^2,\C P^2)) \cong \Z_2\{\nu_5^2 \circ \Sigma^7p\}.
$ \\
$
\pi_8(map_\ast(\C P^2,\C P^2)) \cong \Z_4\{\overline{\nu_5 \circ \eta_8^2}\}\oplus \Z_3\{ {\alpha_2(5)} \circ \Sigma^8 p\} \oplus \Z_5\{ {\alpha_1'(5)} \circ \Sigma^8 p\}
$. \\
$
\pi_9(map_\ast(\C P^2,\C P^2)) \cong \Z_4\{\nu_5 \circ g_8(\C)\}
$. \\
$
\pi_{10}(map_\ast(\C P^2,\C P^2)) \cong \Z_2\{\nu_5^3 \circ \Sigma^{10}p\} \oplus \Z_4\{\overline{\s'''}\}  \oplus \Z_3\{ \overline{\alpha_2(5)} \} \oplus \Z_5\{ \overline {\alpha_1'(5)}  \}
$.\\
$
\pi_{11}(map_\ast(\C P^2,\C P^2)) \cong  \Z_4\{\nu_5 \circ \s_8 \circ \Sigma^{11}p\} \oplus \Z_9\{\beta_1(5) \circ \Sigma^{11}p\}.
$ \\
$
\pi_{12}(map_\ast(\C P^2,\C P^2)) \cong \Z_4\{\zeta_5 \circ \Sigma^{12}p\} \oplus \Z_2^2\{\nu_5 \circ \overline{\nu}_8 \circ \Sigma^{12}p,\nu_5^2 \circ g_{11}(\C)\} \oplus \Z_9\{\alpha_3(5) \circ \Sigma^{12}p\} \oplus \Z_7.
$\\
$
\pi_{13}(map_\ast(\C P^2,\C P^2)) \cong \Z_4\{\nu_5\circ \s_8 \circ \overline{2\iota_{15}}\} \oplus \Z_9\{\overline{\beta_1(5)}\}.
$

\end{theorem}

\section{Applications : Classifying path-components of mapping spaces and Cyclic maps}
In this section, we apply the results obtained in Sections 3, 4 and 5 to classify homotopy types of path-components of certain mapping spaces and compute certain generalized Gottlieb groups.

The term fibration is used for a Hurewicz fibration, that is a map with the homotopy lifting property with respect to all spaces \cite[p.49]{M}.
It is well known that the evaluation map $w_f:map(X,Y;f) \to Y$,
$w_f(g)=g(\ast)$, is a fibration \cite[Lemma 8.15]{G}. For
fibrations $p:E_1 \to B$ and $q:E_2 \to B$, $p$ and $q$ are said
to be {\it fiber homotopy equivalent} if there is a homotopy
equivalence $h:E_1 \to E_2$ such that $q \circ h = p$
\cite[p.52]{M}.

Here we remind several results of the generalized Whitehead product \cite{A}.
If $\alpha$ and $\beta$ are homotopy classes, then the
Whitehead product of $\alpha$ and $\beta$ is denoted by
$[\alpha,\beta]$.
A map $f:X \to Y$ is cyclic if there is a map $F: X \times Y \to Y$, called affiliated map, such that
$F(x,\ast)=f(x)$ and $F(\ast,y)=y$.
We denote the set of cyclic map from $X$ to $Y$ by $G(X,Y)$ and it has group structure if $X$ is a co-H-group \cite{Vara}.
\\

We recall the following equivalent statements due to
\cite[Lemma 2]{Hansen} and \cite[Theorem 3.7]{LS}:

\begin{thm}
Let $\Sigma X$ and $\Sigma Y$ be CW complexes with non-degenerate basepoints and $\Sigma X$ is finite CW complex.
Then the following are equivalent.\\
(A) the map $f:\Sigma X \to \Sigma Y$ is cyclic\\
(B) $[f,id_{\Sigma Y}]=0$ where $[\,\,,\,\,]$ is the generalized Whitehead product\\
(C) the evaluation fibration $w_f:map(\Sigma X,\Sigma Y;f) \to \Sigma Y$ has a section\\
(D) the evaluation fibration $w_f:map(\Sigma X,\Sigma Y;f) \to \Sigma Y$ is fibre-homotopy equivalent to $w_\ast:Map(\Sigma X,\Sigma Y;\ast) \to Y$
\end{thm}

Here is a connection of path components of mapping spaces and cyclic maps.

\begin{thm} \cite[Theorem 3.10]{LS2008}
Suppose $X$ is a CW co-H-space and $Y$ is any CW complex.
Let $d \in G(X,Y)$ be any cyclic map.
Then for each map $f:X\to Y$, we have $map(X,Y;f) \simeq map(X,Y;f+d)$.
If $X$ is a finite co-H-space then the corresponding evaluation fibrations $w_f$ and $w_{f+d}$ are fibre-homotopy equivalent.
\end{thm}

The following Theorem shows a relation between the generaized Whitehead product and evlauation fibration \cite[Theorem 1]{Hansen}

\begin{thm}
Given a pair of homotopy classes $\alpha=[f],\beta=[g] \in [\Sigma A,\Sigma B]$ such that at least one of the identities $[\alpha,\iota_{\Sigma B}] = \pm [\beta,\iota_{\Sigma B}]$ holds.
Then the evaluation fibrations
$w_f:\text{map}(\Sigma A,\Sigma B;f) \to \Sigma B)$
and
$w_g:\text{map}(\Sigma A,\Sigma B;g) \to \Sigma B)$
are fibre homotopy equivalent.
\end{thm}

%
%
%
%
%

We recall Proposition 4.4 of \cite{KMNST}.

\begin{prop}
(1) $[\Sigma \C P^2,S^2] \cong \Z\{\eta_2 \circ \overline{2\iota_3}\}$.\\
(2) $[\Sigma^2 \C P^2,S^3] \cong \Z_2\{\nu' \circ \Sigma^2 p\} \oplus \Z_3\{\alpha_1(3) \circ \Sigma^2 p\}$.\\
(3) $[\Sigma^3 \C P^2,S^4] \cong \Z\{\nu_4 \circ \Sigma^3p\} \oplus \Z_2\{\Sigma \nu' \circ \Sigma^3 p\} \oplus \Z_3\{\alpha_1(4) \circ \Sigma^3p\}$.\\
(4) $[\Sigma^n \C P^2,S^{n+1}] \cong \Z_4\{\nu_{n+1} \circ \Sigma^np\} \oplus \Z_3\{\alpha_1({n+1}) \circ \Sigma^n p\}$ for $n \ge 4$.
\end{prop}

Let $G$ be an abelian group and let $S$ be a subset of $G$.
We denote the smallest subgroup of $G$ containing a subset $S \subset G$ by $\left< S \right >$.
We denote $\omega_f:map(X,Y;f) \to Y$ by $\omega_f$ for $f:X \to Y$.

\begin{thm}
(1) For each $[f],[g]\in [\Sigma \C P^2,S^2]$, the evaluation fibrations
$w_f $ and $w_g $
are fibre homotopy equivalent.\\
(2) For each $[f],[g]\in [\Sigma^2 \C P^2,S^3]$, the evaluation fibrations
$w_f $ and $w_g $
are fibre homotopy equivalent.\\
(3) 
(a) For each $[f],[g] \in \left< 2\nu_4 \circ \Sigma^3p,\Sigma \nu' \circ \Sigma^3p \right>$,
the evaluation fibrations
$w_f $ and $w_g $
are fibre homotopy equivalent.\\
(b) For each $[f] \in \left< 2\nu_4 \circ \Sigma^3p,\Sigma \nu' \circ \Sigma^3p \right>$ and
$[g]  \in \{(2n+1)\nu_4 \circ \Sigma^3p\,|\,n\in \Z\}$,
the evaluation fibrations
$w_g $ and
$w_{f+g} $
are fibre homotopy equivalent.\\
(c) For each $[f]  \in  \left< 2\nu_4 \circ \Sigma^3p,\Sigma \nu' \circ \Sigma^3p \right>$ and
$[h] \in \{ \alpha_1(4) \circ \Sigma^3 p, 2\alpha_1(4) \circ \Sigma^3 p\}$,
the evaluation fibrations
$w_h $ and $w_{h+f} $
are fibre homotopy equivalent.\\
(4) For each $[f],[g]\in [\Sigma^4 \C P^2,S^5]$, the evaluation fibrations\\
$w_f $ and $w_g $
are fibre homotopy equivalent.\\
(5) 
(a) The evaluation fibrations
$w_\ast $ and
$w_f $
are fiber homotopy equivalent where $[f]=2\nu_6 \circ \Sigma^5 p$.\\
(b) For each $[f] \in \{ \nu_6 \circ \Sigma^5 p,3 \nu_6 \circ \Sigma^5 p\}$ and $[g] \in \left<2\nu_6 \circ \Sigma^5 p \right>$,
the evaluation fibrations
$w_f $ and $w_{f+g} $
are fibre homotopy equivalent.\\
(c)  For each $[f] \in \{ \alpha_1(6) \circ \Sigma^5 p,2\alpha_1(6) \circ \Sigma^5 p\}$ and $[g] \in \left<2\nu_6 \circ \Sigma^5 p \right>$,
the evaluation fibrations
$w_f $ and $w_{f+g} $
are fibre homotopy equivalent.\\
(6) For each $[f],[g]\in [\Sigma^6 \C P^2,S^7]$, the evaluation fibrations
$w_f $ and $w_g $
are fibre homotopy equivalent.\\
(7) 
(a) Two evaluation fibrations
$w_f $
and  $w_g $
are fibre homotopy equivalent
where $[f]=\nu_8 \circ \Sigma^7p$ and $[g]=3\nu_8 \circ \Sigma^7p$.\\
(b) Two evaluation fibrations
$w_f $
and
$w_g $
are fibre homotopy equivalent
where $[f]=\alpha_1(9) \circ \Sigma^7p$ and $[g]=2\alpha_1(9) \circ \Sigma^7p$.\\
(8) 
(a) For each $[f],[g] \in \left< \alpha_1(9) \circ \Sigma^8p,2\nu_9 \circ \Sigma^8p \right>$,
the evaluation fibrations
$w_f $
and $w_g $
are fibre homotopy equivalent.\\
(b) For each $[f] \in \{\nu_9 \circ \Sigma^8p,3\nu_9 \circ \Sigma^8p\}$ and $[g] \in\left< \alpha_1(9) \circ \Sigma^8p,2\nu_9 \circ \Sigma^8p \right>$,
the evaluation fibrations
$w_f $ and   $w_{f+g} $ are fibre homotopy equivalent.\\

\begin{proof}
(1)
By \cite[(13)]{KM} and \cite[(2.9)]{BH}, we have $[\iota_2,\eta_2 \circ \overline{2\iota_3}]=0$.
Then we obtain (1) by Theorem 1.\\
(2)
Since $S^3$ is an H-space,
$
[\iota_3,f]=0
$
for any $f\in [\Sigma^2 \C P^2,S^3].$
Then we obtain (2) by Theorem 1.\\
(3)
By \cite[(5.13)]{T} and \cite[Proposition 3.4]{KMNST}, we have
 \begin{eqnarray}
\notag [\iota_4, \nu_4 \circ \Sigma^3 p]  &=& [\iota_4,\nu_4] \circ \Sigma(\iota_3 \wedge \Sigma^2 p)
= 2 \nu_4^2 \circ \Sigma^6 p   \\
\notag  [\iota_4,\Sigma \nu' \circ \Sigma^3 p] &=& [\iota_4,\Sigma \nu'] \circ \Sigma(\iota_3 \wedge \Sigma^2 p)
= 4 \nu_4^2 \circ \Sigma^6p = 0.
 \end{eqnarray}
Then we obtain (a) and (b) by Theorem 1 and 2.
By \cite[Proposition 3.4]{KMNST}
$$
[\iota_4,\alpha_1(4) \circ \Sigma^3p] = [\iota_4,\alpha_1(4)] \circ \Sigma(\iota_3 \wedge \Sigma^2 p)
= [\iota_4,\iota_4] \circ \alpha_1(7) \circ \Sigma^6 p
$$
 is of order 3.
Then we have (c) by Theorem 2 and 3.\\
%
(4)
By \cite[P.48]{T}, we have
$$
[\iota_5,\nu_5 \circ \Sigma^4 p]
=[\iota_5,\nu_5] \circ \Sigma(\iota_4 \wedge \Sigma^3p)
=0
$$
Since $G_8(S^5)=\pi_8(S^5)$ \cite[p.428]{GM1}, we obtain
$$
[\iota_5, \alpha_1(5) \circ \Sigma^4 p]
=[\iota_5, \alpha_1(5) ] \circ \Sigma(\iota_4 \wedge   \Sigma^3p)=0.
$$
Thus we have (4) by Theorem 2 and 3.\\
(5)
By \cite[p.63]{T}, we have
$$
[\iota_6,\nu_6 \circ \Sigma^5 p] = [\iota_6,\nu_6] \circ \Sigma(\iota_5 \wedge \Sigma^4p)
= 2 \overline{\nu}_6 \circ \Sigma^{10} p  .
$$
Since $\overline{\nu}_6 \circ \Sigma^{10} p$ is of order 4  \cite[Proposition 3.7]{KMNST},
we have (a) and (b) by Theorem 1 and 2.
By \cite[Proposition 3.7]{KMNST},
$$
[\iota_6,\alpha_1(6) \circ \Sigma^5p]=[\iota_6,\alpha_1(6)] \circ \Sigma(\iota_5 \wedge \Sigma^4p)=[\iota_6,\iota_6]\circ \alpha_1(11) \circ \Sigma^{10}p
$$
is of order 3.
Then we have (c) by Theorem 2 and 3.\\
%
(6) Since $S^7$ is an H-space,
$
[\iota_7,f]=0
$
for any $f \in[\Sigma^8 \C P^2,S^7].$
Then we obtain (6) by Theorem 1.\\
(7)
By \cite[(7.19)]{T}, we have
$$
[\iota_8,\nu_8 \circ \Sigma^7p] = [\iota_8,\nu_8] \circ \Sigma(\iota_7 \wedge \Sigma^6p)
=2\s_8 \circ \nu_{15} \circ \Sigma^{14}p -x \nu_8 \circ \s_{11} \circ \Sigma^{14}p
$$
for an odd $x$.
By Theorem 1 and 2, we have (a) and (b).
By Proposition 2,
$$
[\iota_8,\alpha_1(8) \circ \Sigma^7 p] = [\iota_8,\alpha_1(8)] \circ \Sigma(\iota_7 \wedge \Sigma^6p) = [\iota_8,\iota_8] \circ \alpha_1(15)\circ \Sigma^{14}p
$$
is of order 3.
By Theorem 2 and 3, we have (c).\\
(8)
By \cite[(7.22)]{T}, we have
$$[\iota_9,\nu_9 \circ \Sigma^8 p] = [\iota_9,\nu_9] \circ \Sigma(\iota_8 \wedge \Sigma^7p)
=\overline{\nu}_9 \circ \nu_{17} \circ \Sigma^{16}p . $$
Since $\overline{\nu}_9 \circ \nu_{17} \circ \Sigma^{16}p$ is of order 2, we have
$[\iota_9,2\nu_9 \circ \Sigma^8 p]= 0 $.
Since $G_{12}(S^9)=\pi_{12}(S^9)$ \cite[p.428]{GM1}, we obtain
$$
[\iota_9,\alpha_1(9)\circ \Sigma^8p] = [\iota_9,\alpha_1(9)] \circ \Sigma(\iota_8 \wedge \Sigma^7p)=0.
$$
By Theorem 2 and 3, we have (b).

\end{proof}
\end{thm}

\begin{cor}
(1) $G_1(\C P^2,S^2) = Z\{\eta_2 \circ \overline{2\iota_3}\}$.\\
(2) $G_2( \C P^2,S^3) = \Z_2\{\nu' \circ \Sigma^2 p\} \oplus \Z_3\{\alpha_1(3) \circ \Sigma^2 p\}$.\\
(3) $G_3( \C P^2,S^4) = \Z\{2\nu_4 \circ \Sigma^3p\} \oplus \Z_2\{\Sigma \nu' \circ \Sigma^3p\}$.\\
(4) $G_4( \C P^2,S^5) = \Z_4\{\nu_{5} \circ \Sigma^4p\} \oplus \Z_3\{\alpha_1(5) \circ \Sigma^4 p\}$. \\
(5) $G_5( \C P^2,S^6) = \Z_2\{2\nu_{6} \circ \Sigma^5p\} $.\\
(6) $G_6( \C P^2,S^7) = \Z_4\{\nu_7 \circ \Sigma^6p\} \oplus \Z_3\{\alpha_1(7) \circ \Sigma^6 p\}$. \\
(5) $G_7( \C P^2,S^8) =0$.\\
(6) $G_8(\C P^2,S^9) = \Z_2\{2\nu_9 \circ \Sigma^8p\} \oplus \Z_3\{\alpha_1(9) \circ \Sigma^8 p\}$.
\end{cor}

\begin{thm}\cite[Corollary 2.4]{LS}
For each $n \geq 1$. we have an isomorphism of abelian groups
$$G_n(map(X,Y;0)) \cong G_n(Y) \oplus G(\Sigma^n X,Y).$$
\end{thm}

By \cite[Corollary 2.4]{LS},
for each $ n,m \geq 1$, we have an isomorphism of abelian groups
$$G_n(map(\Sigma^{m} \C P^2,S^{m+1};0)) \cong G_n(S^{m+1}) \oplus G_m(\C P^2,S^{m+1}).$$

Mukai-Golasinski \cite{GM1} determined Gottelieb groups of sphere $G_{n+k}(S^n)$ for $1 \le k \le 13$, $2 \le n \le 26$.
Therefore we have group structures of
$$G_n(map(\Sigma^m \C P^2,S^{m+1}:0)) \cong G_n(S^{m+1}) \oplus G_m(\C P^2,S^{m+1})$$
for $1 \le m \le 8$ and $1 \le m-n+1 \le 13$.
\\

We recall a long exact sequence \cite{Lang}.
\begin{thm}
There is a long exact sequence
$$\hspace{-20mm}\cdots \to [\Sigma^r A,map(\Sigma B,X;f)] \xrightarrow{w_\ast} [\Sigma^r A,X ] \xrightarrow{P_\alpha}$$
$$\hspace{7mm}\xrightarrow{P_\alpha}[\Sigma(\Sigma^{r-1}A \wedge B),X] \xrightarrow{i_\ast'} [\Sigma^{r-1} A,map(\Sigma B,X;f)]\to \cdots$$
where $\alpha=[f] \in [\Sigma B,X]$,
$P_\alpha(\beta)=[\beta,\alpha]$ is the generalized Whitehead product
and $w:map(\Sigma B,X;f) \to B$ is the evaluation map.\\
\end{thm}

From Theorem 4 and 6, we have the following.

\begin{thm}
(1) $\pi_1(\text{map}(\Sigma \C P^2,S^2;f) )
\cong \Z_2\{\eta_2 \circ \nu' \circ \Sigma^2p\} \oplus \Z_3\{[\iota_2,\iota_2] \circ \alpha_1(3)\}$
for all $f : \Sigma \C P^2 \to S^2$.\\

(2)  $\pi_2(\text{map}(\Sigma^2 \C P^2,S^3;f))
\cong \Z_2\{\nu' \circ \overline{2\iota_6}\} \oplus \Z_3\{\overline{\alpha_1(3)}\}$
for all $f : \Sigma^2 \C P^2 \to S^3$.\\

(3)
(a) For $[f] \in \left< 2\nu_4 \circ \Sigma^3p,\Sigma \nu' \circ \Sigma^3p \right>$,
$\pi_3(\text{map}(\Sigma^3 \C P^2,S^4;f))
\cong \Z_4\{\nu_4^2 \circ \Sigma^6 p\} \oplus \Z_3^2\{\alpha_1(4) \circ \alpha_1(7) \circ \Sigma^6p,[\iota_4,\iota_4] \circ \alpha_1(7) \circ \Sigma^6p\}$.\\

(b) For $[f] \in \{(2n+1)\nu_4 \circ \Sigma^3p\,| \, n \in \Z\}$
and
$[g] \in \left< 2\nu_4 \circ \Sigma^3p,\Sigma \nu' \circ \Sigma^3p \right>$,
$\pi_3(\text{map}(\Sigma^3 \C P^2,S^4;f+g))
\cong \Z_2\{\nu_4^2 \circ \Sigma^6 p\} \oplus \Z_3^2\{\alpha_1(4) \circ \alpha_1(7) \circ \Sigma^6p,[\iota_4,\iota_4] \circ \alpha_1(7) \circ \Sigma^6p\}$.\\

(c) For $[f] \in \{\alpha_1(4) \circ \Sigma^3p, 3\alpha_1(4) \circ \Sigma^3p\}$
and
$[g] \in \left< 2\nu_4 \circ \Sigma^3p,\Sigma \nu' \circ \Sigma^3p \right>$,
$\pi_3(\text{map}(\Sigma^3 \C P^2,S^4;f+g))
\cong \Z_4\{\nu_4^2 \circ \Sigma^6 p\} \oplus \Z_3\{\alpha_1(4) \circ \alpha_1(7) \circ \Sigma^6p\}$.\\

(4) $\pi_4(\text{map}(\Sigma^4 \C P^2,S^5;f))
\cong \Z_4\{\overline{\nu_5 \circ \eta_8^2} \} \oplus \Z_3\{\alpha_2(5) \circ \Sigma^8 p\} \oplus \Z_5\{\alpha_{1,5}(5) \circ \Sigma^8p\}$
for all $f : \Sigma^4 \C P^2 \to S^5$.\\

(5)
(a) For $[f] \in \left< 2\nu_6 \circ \Sigma^5p \right>$,
$\pi_5(\text{map}(\Sigma^5 \C P^2,S^6;f))
\cong \Z_4^2\{\nu_6 \circ g_9, \overline{\nu}_6 \circ \Sigma^{10}p\} \oplus \Z_3\{[\iota_6,\iota_6] \circ \alpha_1(11) \circ \Sigma^{10}p\}$.\\

(b) For $[f] \in \{\nu_6 \circ \Sigma^5p, 3\nu_6 \circ \Sigma^5p\}$ and $[g] \in \left< 2\nu_6 \circ \Sigma^5p \right>$,
$\pi_5(\text{map}(\Sigma^5 \C P^2,S^6;f))$\\
$\cong \Z_4\{\nu_6 \circ g_9\} \oplus \Z_2\{ \overline{\nu}_6 \circ \Sigma^{10}p\} \oplus \Z_3\{[\iota_6,\iota_6] \circ \alpha_1(11) \circ \Sigma^{10}p\}$.\\

(c) For $[f] \in \{\alpha_1(6) \circ \Sigma^5p, 2\alpha_1(6) \circ \Sigma^5p\}$
and
$[g] \in \left< 2\nu_6 \circ \Sigma^5p \right>$,\\
$\pi_5(\text{map}(\Sigma^5 \C P^2,S^6;f))
\cong \Z_4^2\{\nu_6 \circ g_9, \overline{\nu}_6 \circ \Sigma^{10}p\}$.\\

(6) $\pi_6(\text{map}(\Sigma^6 \C P^2,S^7;f))
\cong \Z_8\{\s' \circ \overline{2\iota_{14}}\} \oplus \Z_3\{\overline{\alpha_2(7)}\} \oplus \Z_5\{\overline{\alpha_{1,5}(7)}\}$
for all $f : \Sigma^6p \to S^7$.\\

(7)
(a) $\pi_7(\text{map}(\Sigma^7 \C P^2,S^8;\ast))
\cong \Z_4^2\{\s_8 \circ \nu_{15} \circ \Sigma^{14}p,\nu_8 \circ \s_{11} \circ \Sigma^{14}p\}
\oplus \Z_3^2\{\beta_1(8) \circ \Sigma^{14}p,[\iota_8,\iota_8]\circ \alpha_1(15) \circ \Sigma^{14}p\}$.\\

(b) For $[f] \in \{\nu_8 \circ \Sigma^7p, 3\nu_8 \circ \Sigma^7p\}$,
$\pi_7(map(\Sigma^7 \C P^2, S^8;f))\cong
\Z_4\{\s_8 \circ \nu_{15} \circ \Sigma^{14}p +\nu_8 \circ \s_{11} \circ \Sigma^{14}p\} \oplus \Z_3^2\{[\iota_8,\iota_8]\circ \alpha_1(5) \circ \Sigma^{14}p,\beta_1(8) \circ \Sigma^{14}p\}.$\\

(c) For $[f] = 2\nu_8 \circ \Sigma^7 p$,
$\pi_7(\text{map}(\Sigma^7 \C P^2,S^8;f))
\cong \Z_4\{\s_8 \circ \nu_{15} \circ \Sigma^{14}p\} \oplus \Z_2\{\nu_8 \circ \s_{11} \circ \Sigma^{14}p\}
\oplus \Z_3^2\{\beta_1(8) \circ \Sigma^{14}p,[\iota_8,\iota_8]\circ \alpha_1(15) \circ \Sigma^{14}p\}$.\\

(d) $[f] \in \{\alpha_1(9) \circ \Sigma^7 p, 2\alpha_1(9) \circ \Sigma^7 p\}$,
$\pi_7(\text{map}(\Sigma^7 \C P^2,S^8;f))
\cong \Z_4^2\{\s_8 \circ \nu_{15} \circ \Sigma^{14}p,\nu_8 \circ \s_{11} \circ \Sigma^{14}p\}
\oplus \Z_3\{\beta_1(8) \circ \Sigma^{14}p\}$.\\

(8)
(a) For $[f] \in \left< \alpha_1(9) \circ \Sigma^8p, 2\nu_9 \circ \Sigma^8p \right>$,
$\pi_8(\text{map}(\Sigma^8 \C P^2,S^9;f))
\cong \Z_8\{\nu_9^2 \circ g_{15}\} \oplus \Z_2^2\{\overline{\nu}_9 \circ \nu_{17} \circ \Sigma^{16}p, \overline{\eta_9 \circ \e_{10}}\}
\oplus \Z_7\{\alpha_{1,7}(9) \circ \Sigma^{14}p\} \oplus \Z_9\{\alpha_3'(9) \circ \Sigma^{14}p\}$.\\

(b) For $[f] \in \{\nu_9 \circ \Sigma^8p,3\nu_9 \circ \Sigma^8p \}$
and
$[g] \in \left< \alpha_1(9) \circ \Sigma^8p, 2\nu_9 \circ \Sigma^8p \right>$,\\
$\pi_8(\text{map}(\Sigma^8 \C P^2,S^9;f+g))
\cong \Z_2\{\nu_9^2 \circ g_{15}\} \oplus \Z_2\{\overline{\eta_9 \circ \e_{10}}\}
\oplus \Z_7\{\alpha_{1,7}(9) \circ \Sigma^{14}p\} \oplus \Z_9\{\alpha_3'(9) \circ \Sigma^{14}p\}$.
\end{thm}

We denote the path-components of mapping space $map(X,Y)$ containing $f$ by $map( f)$. \\

Finally, we classify homotopy types of path-components of the following mapping spaces.

\begin{thm}
(1)
$map(\Sigma^n \C P^2, S^{n+1};f)$ consists of one path-component for $n=1,2,4,6$.\\


(2)
$map(\Sigma^3 \C P^2, S^4;f)$ consists of four path-components
$map(f)$, $map(f+g)$, $map(f+h)$ and $map(f+g+h)$
where
$[f] \in \left< 2\nu_4 \circ \Sigma^3p,\Sigma \nu' \circ \Sigma^3p \right>$,
$[g] \in \{(2n+1)\nu_4 \circ \Sigma^3p\,| \, n \in \Z\}$
and
$[h] \in \{\alpha_1(4) \circ \Sigma^3p, 2\alpha_1(4) \circ \Sigma^3p\}$.\\

(3)
$map(\Sigma^5 \C P^2, S^6;f)$ consists of four path-components
$map(f)$, $map(f+g)$, $map(f+h)$ and $map(f+g+h)$
where
$[f] \in \left< 2\nu_6 \circ \Sigma^5p \right>$,
$[g] \in \{\nu_6 \circ \Sigma^5p, 3\nu_6 \circ \Sigma^5p\}$
and
$[h] \in \{\alpha_1(6) \circ \Sigma^5p, 2\alpha_1(6) \circ \Sigma^5p\}$.\\

(4)
$map(\Sigma^7 \C P^2, S^8;f)$ consists of six path-components
$map(\ast)$, $map(f)$, $map(2f)$, $map(g)$, $map(f+g)$ and $map(2f+g)$
where
$[f] = \nu_8 \circ \Sigma^7 p$ and  $[g] = \alpha_1(8) \circ \Sigma^7 p$.\\

(5)
$map(\Sigma^8 \C P^2, S^9;f)$ consists of two path-components
$map(f)$ and $map(f+g)$
where
$[f] = \nu_8 \circ \Sigma^7 p$ and  $[g] = \alpha_1(8) \circ \Sigma^7 p$.
\end{thm}


\end{document}